\documentclass{amsart}
\usepackage{amscd, amsfonts}
\newtheorem{thm}{Theorem}

\newtheorem{cor}[thm]{Corollary}
\newtheorem{conj}[thm]{Conjecture}

\theoremstyle{remark}
\newtheorem{example}[thm]{Example}
\newtheorem{remark}[thm]{Remark}

\theoremstyle{definition}
\newtheorem{definition}[thm]{Definition}

\newcommand{\into}{\hookrightarrow}

\DeclareMathOperator*{\colim}{colim}

\begin{document} 
\title[A duality between string topology and fusion]{A duality between string topology and the fusion product in equivariant K-theory}
\author{Kate Gruher}
\address{Stanford University}
\email{gruher@math.stanford.edu}
\thanks{The author was supported by a National Defense Science and Engineering Graduate Research Fellowship}
\date{\today}

\bibliographystyle{amsalpha}
\begin{abstract}
Let $G$ be a compact Lie group.  Let $G \to E \to M$ be a principal $G$-bundle over a closed manifold $M$, and $G \to Ad(E) \to M$ its adjoint bundle.  In this paper we describe a new Frobenius algebra structure on $h_*(Ad(E)),$ where $h_*$ is an appropriate generalized homology theory.  Recall that a Frobenius algebra has both a product and a coproduct.  The product in this new Frobenius algebra is induced by the string topology product.  In particular, the product can be defined when $G$ is any topological group and in the case that $E$ is contractible it is precisely the Chas-Sullivan string product on $H_*(LM).$  We will show that the coproduct is induced by the Freed-Hopkins-Teleman fusion product.  Indeed, when $M$ is replaced by $BG$ and $h_*$ is $K$-theory the coproduct is the completion of the Freed-Hopkins-Teleman fusion structure.  We will then show that this duality between the string and fusion products is realized by a Spanier-Whitehead duality between certain Thom spectra of virtual bundles over $Ad(E)$.

\end{abstract}

\maketitle

\section{Introduction}

Let $G$ be a topological group.  Consider a principal $G$-bundle $G \to E \to M$ over a closed manifold $M$ and let $G \to Ad(E) \to M$ be its associated adjoint bundle.  The author and Salvatore \cite{gs} recently decribed a string product on $h_*(Ad(E)),$ analogous to the Chas-Sullivan string product, where $h_*$ is a generalized homology theory that supports an orientation of $M$.  When $E$ is contractible, $Ad(E) \simeq LM$ and this product is the original Chas-Sullivan product on $H_*(LM)$ \cite{cs}.  Using the techniques of Cohen and Jones \cite{cj}, it was shown that the Thom spectrum $Ad(E)^{-TM}$ has a ring spectrum structure that realizes the string product.

Now let $G$ be a compact, connected, simply connected, simple Lie group.  Freed, Hopkins, and Teleman recently described a product in the twisted equivariant $K$-cohomology of $G$, ${}^\tau \! K_G^*(G),$ which realizes the fusion product on the Verlinde algebra of representations of the loop group $LG$ \cite{fht3}, \cite{freedv}.  In this paper we first observe that the FHT construction also defines a product on ${}^\tau \! K^*(Ad(E)),$ for certain twistings $\tau$, and in the case $\tau=0$ the classifying map $f_E: M \to BG$ induces a ring homomorphism $f_E^*: K_G^*(G) \to K^*(Ad(E)).$

We further show that for any compact Lie group $G$, the FHT construction defines a product on $h^*(Ad(E))$ when $h^*$ is an appropriate untwisted generalized cohomology theory.  We refer to this product also as the fusion product.  We prove a duality result between the string product in $h$-homology and the fusion product in $h$-cohomology.  This duality leads to the new Frobenius algebra structure on $h_*(Ad(E)$.  We show that this duality structure arises from a Spanier-Whitehead duality between two Thom spectra of virtual bundles over $Ad(E).$  We now state these results more precisely.

Let $G$ be a compact Lie group and let $h^*$ be a multiplicative generalized cohomology theory.

\begin{definition}
The group $G$ has an $h^*$-adjoint orientation if the vertical tangent bundle of $Ad(EG) \to BG$ is $h^*$-oriented.
\end{definition}

\begin{remark}
If $G$ is connected, then it admits an $H^*$-adjoint orientation.  If in addition $G$ is simply connected, then it admits a $K^*$-adjoint representation.
\end{remark}

Let $G \to E \to M$ be a principal $G$-bundle over a closed manifold $M$.  For the following theorem, let $h^*$ be a multiplicative cohomology theory such that the coefficient ring $h^*(\ast)$ is a graded field.  Let $h_*$ be the associated homology theory. 
\begin{thm}\label{frob}

If $M$ is $h^*$-oriented and $G$ is $h^*$-adjoint oriented, then $h_*(Ad(E))$ has a Frobenius algebra structure over $h_*(\ast)$ with multiplication given by the string product and comultiplication given by the fusion product. 

\end{thm}
Both the CS string product and the FHT fusion product have given rise to interesting Frobenius algebras.  In \cite{cg}, Cohen and Godin describe a unital, non-counital, commutative Frobenius algebra structure on $h_*(LM)$, with the multiplication given by the string product.  In \cite{fht3}, Freed, Hopkins, and Teleman describe a Frobenius algebra structure on ${}^{\tau}\! K_G^*(G)$, for transgressed twistings, that is isomorphic to a Frobenius algebra structure on the Verlinde algebra. Since Frobenius algebras correspond to $2$-dimensional topological field theories, this Frobenius algebra structure on the twisted $K$-theory of $G$ is called the ``Verlinde TFT.''  The Frobenius algebra structure in Theorem~\ref{frob} is different from both of the ones described above.

Our next result realizes the duality between the string and fusion products on the level of spectra.  Recall the following theorem from \cite{gs}.
\begin{thm}\label{genstring}
$Ad(E)^{-TM}$ is a ring spectrum.  If $M$ is oriented, then the induced product on $H_*(Ad(E))$ is the string product.
\end{thm}
The spectrum denoted $Ad(E)^{-TM}$ is the Thom spectrum of the virtual bundle $-p^*(TM)$, where $p:Ad(E) \to M$.  Since $Ad(E)^{-TM}$ is a finite spectrum, it has a Spanier-Whitehead dual which then has a coproduct.  Denote by $T_v$ the vertical tangent bundle of $Ad(E) \to M$.  We can form the bundle $ad(E) = E \times _G \mathfrak{g}$ over $M$, where the Lie algebra $\mathfrak{g}$ is given the adjoint action; then $T_v \cong p^*(ad(E)).$  Notice also that $T_v$ is the pullback of the vertical tangent bundle of $Ad(EG) \to BG$ via the classifying map of $E$.   The spectrum-level duality is described by the following theorem.

\begin{thm}\label{dual}

The Spanier-Whitehead dual of $Ad(E)^{-TM}$ is $Ad(E)^{-T_v}.$  When $G$ has an $h^*$-adjoint orientation, the product on $h^*(Ad(E))$ given by applying Thom isomorphisms to the induced co-ring spectrum structure on $Ad(E)^{-T_v}$ is the fusion product.

\end{thm}

The string topology of $BG$ was defined in \cite{gs} as a pro-ring spectrum $LBG^{-TBG}$.  The ring spectra appearing in the construction of $LBG^{-TBG}$ are all of the form $Ad(E)^{-TM}.$  Theorem~\ref{dual} describes a duality between the string topology of $BG$ and the fusion product on $h^*(Ad(EG)).$  This setting is discussed more in section~\ref{section:further}.

Before proceeding, let us consider two examples.  Theorem~\ref{frob} interpolates between the Frobenius algebras given by the two extreme cases of the theorem: $G = \{ e \}$ and $M = \ast$.  Let $k$ be a field.

%
\begin{example} $G = \{ e \}.$
Then $E = Ad(E) = M$.  The multiplication on $H_*(M; k)$ is given by the intersection product and the comultiplication by the diagonal $\Delta _*$.  The nondegenerate bilinear form for the Frobenius algebra
$$\beta: H_*(M;k) \otimes H_*(M;k) \to k$$
 is given by 
$$x\otimes y \mapsto <Dx,y>$$
 where $Dx$ is the Poincare dual of $x$.
\end{example}

\begin{example} $M = \ast.$
Then $E = Ad(E) = G$.
The multiplication on $H_*(G; k)$ is given by $m_*$ and the comultiplication by $m_!$.  In this case $\beta$ is given by
$$\beta(x\otimes y) = \left\{ \begin{array}{c@{\quad if \quad}l}
0 & |x|+|y| \neq \mbox{dim}(G) \\
m_*(x\times y) &  |x|+|y| = \mbox{dim}(G)
\end{array} \right.
$$
where $H_{\mbox{dim}(G)}(G;k) \cong k$ via $[G] \mapsto 1$.
\end{example}

The relationship between the Frobenius algebra in Theorem~\ref{frob} and the Verlinde TFT will be studied further in \cite{gw}.

This paper will be organized as follows.  In section~\ref{section:cs} we will recall the Chas-Sullivan string product and the Cohen-Jones construction and describe their generalizations to adjoint bundles, including the construction of $LBG^{-TBG}$.  In section~\ref{section:fht} we will recall the main theorem of \cite{fht3} and describe the generalization of the fusion product to adjoint bundles.  In section~\ref{section:compare} we will prove Theorem~\ref{frob} and Theorem~\ref{dual} and give a more detailed description of the Frobenius algebra structure.  In section~\ref{section:further} we will discuss further implications of Theorem~\ref{dual} using Spanier-Whitehead duality for pro-spectra.

The author is grateful to R. Cohen for his guidance and help throughout this project.

\section{String Topology in its Natural Habitat and on Adjoint Bundles}\label{section:cs}

The Chas-Sullivan string product is a product on the homology of the free loop space of a closed, oriented, $d$-dimensional manifold $M$:
$$\mu : H_r(LM) \otimes H_s(LM) \to H_{r+s-d}(LM) $$
that is a combination of the intersection product on the manifold and the Pontrjagin product on the based loop space $\Omega M$.  It can be described as the following composition:
$$H_r(LM) \otimes H_s (LM) \stackrel{\times}{\to} H_{r+s}(LM \times LM) \stackrel{\tilde{\Delta}_!}{\to} H_{r+s-d}(LM \times _M LM) \stackrel{\gamma _*}{\to} H_{r+s-d}(LM) $$
where 
$$LM \times _M LM = \{(\alpha, \beta) \in LM\times LM | \alpha(0) = \beta(0) \},$$
 $\tilde{\Delta}$ is the inclusion $LM\times _M LM \into LM \times LM$, and $\gamma: LM \times _M LM \to LM$ is the concatenation of loops with the same starting point.  

The map $\tilde{\Delta}_!$ is defined using a Pontrjagin-Thom collapse map 
$$\tau_{\tilde{\Delta}}: LM \times LM \to (LM \times _M LM)^{TM}$$
 and the Thom isomorphism.  The construction of $\tau_{\tilde{\Delta}}$ relies on the pullback square of fiber bundles
$$\begin{CD}
LM \times _M LM @>{\tilde{\Delta}}>> LM \times LM \\
@VVV @VVV \\
M @>{\Delta}>> M \times M 
\end{CD} $$
and on the compactness of $M$.  Indeed the CS string product may be generalized to any fiber bundle $W \to M$ where the base $M$ is a compact, oriented manifold and the bundle has the structure of a fiberwise monoid (a monoid in the category of spaces over $M$).  This was done in \cite{gs}.  We refer to resulting product on $H_*(W)$ as the string product.  

Cohen and Jones showed in \cite{cj} that there is a ring spectrum structure on $LM^{-TM}$ that realizes the CS string product in homology.  The multiplication on $LM^{-TM}$ is given by: 
$$LM^{-TM} \wedge LM^{-TM} \simeq (LM\times LM)^{-(TM\times TM)} \stackrel{\tau_{\tilde{\Delta}}}{\to} (LM \times _M LM)^{-TM} \stackrel{\gamma}{\to} LM^{-TM} $$
where $\tau_{\tilde{\Delta}}$ is a Pontrjagin-Thom collapse map, twisted with respect to the virtual bundle $-(TM\times TM)$.  As shown in \cite{gs}, the Cohen-Jones construction also generalizes to fiberwise monoids over compact manifolds, yielding a ring spectrum structure on $W^{-TM}$ that realizes the string product in homology.  

Adjoint bundles of principal bundles form a nice class of fiberwise monoids.  If $G \to E \to B$ is a principal $G$-bundle, then $Ad(E) = E\times _G G = (E \times G) / G$, where $G$ acts on itself by conjugation, is a fiberwise monoid over $B$ (in fact, it is a fiberwise group).  For any topological group $G$, there is a homotopy equivalence of fiberwise monoids over $BG$ between $LBG$ and $Ad(EG)$.  This was used in \cite{gs} to define a pro-ring spectrum $LBG^{-TBG}$ for any compact Lie group $G$, extending the idea of the CS string product to classifying spaces of Lie groups.  $LBG^{-TBG}$ is defined as 
$$ Ad(EG_1)^{-TM_1} \leftarrow \ldots \leftarrow Ad(EG_n)^{-TM_n} \leftarrow Ad(EG_{n+1})^{-TM_{n+1}} \leftarrow \ldots $$
where the $M_n$ are manifolds approximating $BG$ (in particular, $BG = \colim M_n$), and $EG_n$ is the pullback of $EG$ to $M_n$.  The ring spectrum structure on each $Ad(EG_n)^{-TM_n}$ comes from the generalization of the Cohen-Jones construction.

\section{Fusion in its Natural Habitat and on Adjoint Bundles}\label{section:fht}

We will now discuss the previously mentioned result of Freed, Hopkins, and Teleman, and the extension of their construction to the twisted $K$-theory of $Ad(E).$  For this discussion, we will assume that $G$ is a compact, connected, simply connected, simple Lie group.  The simplest example is $G = SU(2).$

 For any $k \in \mathbb{Z}$, let $V_k(G)$ denote the Verlinde algebra of $G$ at level $k$, and let $\zeta(k) = k + h(G) \in H^3_G(G) \cong \mathbb{Z}$, where $h(G)$ is the dual Coxeter number of $G$.  In \cite{fht3}, Freed, Hopkins, and Teleman describe a ring isomorphism
$${}^{\zeta(k)} \! K_G^{|G|}(G) \cong V_k(G)$$
when $k + h(G) > 0$.  The left side is the twisted $K$-theory of $G$, equivariant with respect to the conjugation action, with twisting $\zeta(k)$.  $|G|$ denotes the dimension of $G$ modulo $2$.  See \cite{astwk} for the definition and basic properties of twisted $K$-theory.

Take $\tau$ to be representative cocyle (in some model of equivariant cohomology) for $[\tau] \in H^3_G(G; \mathbb{Z})$.
The ring multiplication $\rho$  on $^\tau \! K_G^*(G)$ is given by the composition
$$^\tau \! K_G^*(G) \otimes {}^\tau \! K_G^*(G) \stackrel{\times}{\to} {}^{(\tau, \tau)} \! K^*_{G\times G}(G\times G) \stackrel{r}{\to} {}^{(\tau, \tau)} \! K^*_{G}(G\times G) \stackrel{m^!}{\to} {}^\tau \! K^*_G(G) $$
where $r$ is the restriction map and $m: G\times G \to G$ is the group multiplication.

The map $m^!$ is an umkehr map for the principal $G$-bundle $G\times G \stackrel{m}{\to} G$. For clarity, let $X=G\times G$ denote the total space of this bundle.   The umkehr map is constructed by using an equivariant Pontrjagin-Thom collapse map $S^V \to G^{\nu}$, where $V$ is a right representation of $G$, to construct an equivariant map $X \times _G S^V \to X \times_G G^\nu$ \cite{bg}.   This map is equivariant with respect to the left action of $G$ given by $g[x,h] = [gxg^{-1}, hg^{-1}]$.  It induces an equivariant map of Thom spaces
$$G^{\eta} \to X^{\beta}$$
where $\eta: X\times _G V \to G$ and $\beta: X\times _G \nu \to X$.  Since $ \beta \oplus T_v = m^*(\eta)$, we have an induced map
$$m^!: {}^{(\tau, \tau)}\!K_G^*(G\times G) \to {}^{\tau}\!K_G^*(G),$$
which shifts the degree by $|G|$.  The shift in the twisting arising from the equivariant $K$-theory Thom isomorphism is zero by our assumptions on $G$.  Notice that since $m^!$ shifts the degree by $|G|$, ${}^{\tau}K_G^{|G|}(G)$ is a subring of ${}^{\tau}K_G^*(G)$, but ${}^{\tau}K_G^{|G|+1}(G)$ is not.

We can extend $\rho$ to the twisted $K$-theory of total spaces of adjoint bundles by defining, for any principal $G$-bundle $E \to B$ and any $\tau \in H^3(Ad(E))$ that is pulled back from $H^3(Ad(EG))$, a product on ${}^{\tau}\!K^*(Ad(E))$:
\begin{multline}\label{rhodef}\rho: {}^{\tau}\!K^*(G \times _G E) \times {}^{\tau}\!K^*(G \times _G E) \stackrel{\times}{\to} {}^{(\tau, \tau)}\!K^*((G \times G)\times  _{G \times G} (E \times E)) \\
 \stackrel{\tilde{\Delta}^*}{\to} {}^{(\tau, \tau)}\!K^*((G \times G) \times _G E) \stackrel{m^!}{\to} {}^{\tau}\!K^*(G \times _G E).
\end{multline}
The construction of $m^!$ in this situation is a non-equivariant version of the Becker-Gottlieb transfer \cite{bg}, yielding a map of Thom spaces
\begin{equation}\label{BGtransfer}
Ad(E)^{\eta} \to (Ad(E) \times _B Ad(E))^{\beta}
\end{equation}
where $\beta \oplus T_v = m^*(\eta)$; here $T_v$ denotes the vertical tangent bundle of 
$$m: Ad(E) \times _B Ad(E) \to Ad(E).$$
By our assumptions on $G$ the bundle $T_v$ is spin and so applying Thom isomorphisms to both sides induces
$$m^!: {}^{(\tau,\tau)}\! K^*((G \times G) \times _G E) \to {}^{\tau}\!K^*(G \times _G E)$$
for any $\tau \in H^3(Ad(E))$ that is pulled back from $H^3(Ad(EG))$.  Again, $m^!$ shifts the degree by $|G|$.
Notice that the construction of $m^!$ relies on the compactness of the group $G$ but does not require the base space of the bundle to be compact.

In the case $\tau =0$, by the Atiyah-Segal completion theorem \cite{as},
 $$K^*(Ad(EG)) \cong K^*_G(G) \hat{}\, ,$$
 the completion of $K^*_G(G)$ at the augmentation ideal.
Then all of the maps
\begin{equation}\label{completion}
K_G^*(G) \to K_G^*(G)\hat{}\, \cong K^*(Ad(EG)) \to K^*(Ad(EG_n))
\end{equation}
are ring maps.

For any compact Lie group $G$, the construction above defines a product on $h^*(Ad(E))$ for any (untwisted) generalized cohomology theory $h^*$ for which $G$ has an $h^*$-adjoint orientation.  Indeed, the map of Thom spaces (\ref{BGtransfer}) induces a map
$$m^!: h^*(Ad(E)\times_B Ad(E)) \to h^*(Ad(E))$$
and the product $\rho$ is defined as in (\ref{rhodef}).

\section{Comparison of the string and fusion products}\label{section:compare}

The construction of the string product on the total space of a fiberwise monoid relies on the compactness of the base.  In particular, the string product can be constructed on adjoint bundles of principal bundles with compact base spaces.  On the other hand, the construction of the fusion product relies on the compactness of the fiber (i.e. the group $G$).  It can be constructed on adjoint bundles of principal $G$-bundles when $G$ is a compact Lie group.

These observations suggest that a natural place to compare the string and fusion products is on adjoint bundles where both the fiber and the base are compact.  If $G\to E \to M$ is a principal bundle where $M$ is a closed oriented manifold and $G$ is a compact Lie group, then we can construct both the string product on the homology of $Ad(E)$ and the fusion product on the cohomology of $Ad(E).$  Notice that this is not an artificial situation: the adjoint bundles $Ad(EG_n)$ appearing in the construction of $LBG^{-TBG}$  are precisely of this form.  Freed, Hopkins, and Teleman were also studying $Ad(EG)$, and we will see that our comparison of the string and fusion constructions on compact adjoint bundles gives rise to a comparison of the constructions on $Ad(EG).$  We now describe this relationship by giving a proof of Theorem~\ref{dual} as stated in the introduction.

\begin{proof}[Proof of Theorem~\ref{dual}]
By \cite{a}, the Spanier-Whitehead dual of $Ad(E)^{-TM}$ is 
$$Ad(E)^{TM - TAd(E)} = Ad(E)^{-T_v}. $$

If $f: M {\to} N$ is a smooth map of manifolds then the Spanier-Whitehead dual of $f$ is the generalized Pontrjagin-Thom map $N^{-TN} \to M^{-TM}$, which is the stable collapse map $N \to M^{\nu}$ twisted with respect to the virtual bundle $-TN$.  The ring spectrum structure on $Ad(E)^{-TM}$ is given by
$$\mu: Ad(E)^{-TM} \wedge Ad(E)^{-TM} \stackrel{\tau_{\tilde{\Delta}}}{\to} (Ad(E) \times _M Ad(E))^{-TM} \stackrel{m}{\to} Ad(E)^{-TM}$$
where $\tau_{\tilde{\Delta}}$ is the Pontrjagin-Thom collapse map for the embedding
$$\tilde{\Delta}: Ad(E)\times _M Ad(E) \into Ad(E) \times Ad(E)$$
twisted with respect to the virtual bundle $-(TM \times TM)$.  The Spanier-Whitehead dual of $\tau_{\tilde{\Delta}}$ is then 
$$ \tilde{\Delta}: (Ad(E)\times _M Ad(E))^{-(T_v \times T_v)} \to (Ad(E) \times Ad(E))^{-(T_v \times T_v)}.$$
The Spanier-Whitehead dual of $m$ is the Pontrjagin-Thom map
$$\tau_m : Ad(E)^{-T_v} \to (Ad(E) \times _M Ad(E))^{-(T_v \times T_v)}.$$
Thus the coproduct on $Ad(E)^{-T_v}$ is given by:
\begin{equation}\label{coprod}
\rho: Ad(E)^{-T_v} \stackrel{\tau_m}{\to} (Ad(E) \times _M Ad(E))^{-(T_v \times T_v)} \stackrel{\tilde{\Delta}}{\to} Ad(E)^{-T_v} \wedge Ad(E)^{-T_v} .
\end{equation}
Applying a multiplicative cohomology theory $h^*$ such that $G$ has an $h^*$-adjoint orientation and then applying the Thom isomorphism on both sides yields the fusion product on $h^*(Ad(E))$ as described in section~\ref{section:fht}.

\end{proof}

We now describe the application of this Spanier-Whitehead duality result to the description of the Frobenius algebra structure where the string product and the fusion product are dual to each other.
Let $k$ be a field.  Recall that a Frobenius algebra is a finite dimensional $k$-algebra $A$ equipped with a nondegenerate bilinear pairing $\beta: A \otimes A \to A$ that is associative, meaning that $\beta(xa\otimes y) = \beta(x\otimes ay)$ for all $x, y, a \in A$.

An equivalent definition is \cite{abrams}:
\begin{definition}
A Frobenius algebra is a $k$-vector space $A$ along with
\begin{itemize}
\item{a multiplication $\mu:A \otimes A \to A$ with unit $\eta: k \to A$; and}
\item{a comultiplication $\rho: A \to A \otimes A$ with counit $\epsilon: A \to k$}
\end{itemize}
such that the Frobenius relation holds:
$$(1 \otimes \mu) \circ (\rho \otimes 1) = \rho \circ \mu = (\mu \otimes 1) \circ (1 \otimes \rho).$$
\end{definition}
The bilinear pairing from the first definition can be written in terms of the maps in the second definition as $\beta = \epsilon \circ \mu$.

Before proving Theorem~\ref{frob} we will give an explicit description of the Frobenius algebra structure on $H_*(Ad(E);k)$:
\begin{itemize}
\item{$\mu = m_* \tilde{\Delta}_!$ (the string product)}
\item{$\rho = \tilde{\Delta}_* m_!$ (the fusion coproduct)}
\item{$\eta: 1 \mapsto s_*[M]$, where $s: M \to Ad(E)$ is the fiberwise monoid unit}
\item{$\epsilon = (pt)_* s_!$, where $pt: M \to \ast$ is the unique map from $M$ to a point and we identify $H_*(\ast ;k) \cong k$.}
\item{$\beta = (pt)_* s_! m_* \tilde{\Delta}_!$}
\end{itemize}

\begin{proof}[Proof of Theorem~\ref{frob}]
Let $h_*$ be a homology theory as in the statement of the theorem.
Combined with Thom isomorphisms, the ring spectrum structure on $Ad(E)^{-TM}$ and its unit $S \to Ad(E)^{-TM}$ define the associative string product $\mu$ on $h_*(Ad(E))$ and unit $\eta: k \to h_*(Ad(E)).$  Likewise, the induced co-ring spectrum structure on $Ad(E)^{-T_v}$ and it counit $Ad(E)^{-T_v} \to S$ define the coassociative fusion coproduct $\rho$ on $h_*(Ad(E))$ and counit $\epsilon: h_*(Ad(E)) \to k.$  We need to check that the Frobenius relation holds, which is equivalent to showing that $\rho$ is a map of $(h_*(Ad(E)), h_*(Ad(E)))$-bimodules.  

By Spanier-Whitehead duality, $Ad(E)^{-T_v} \simeq F(Ad(E)^{-TM}, S)$ has a natural $(Ad(E)^{-TM}, Ad(E)^{-TM})$-bimodule structure up to homotopy.  This module structure coincides with the composition of comultiplication with evaluation; for instance, the right module structure is given by:
\begin{multline*}
F(Ad(E)^{-TM}, S)\wedge Ad(E)^{-TM} \stackrel{\rho \wedge 1}{\to} F(Ad(E)^{-TM} \wedge Ad(E)^{-TM}, S)\wedge Ad(E)^{-TM} \\
 \simeq F(Ad(E)^{-TM}, S)\wedge F(Ad(E)^{-TM}, S) \wedge Ad(E)^{-TM} \stackrel{1 \wedge \theta}{\to} F(Ad(E)^{-TM}, S) \wedge S.
\end{multline*}
Here $\theta: F(Ad(E)^{-TM}, S) \wedge Ad(E)^{-TM} \to S$ is the Spanier-Whitehead duality map.
Applying homology and Thom isomorpisms to this right module structure defines a product $\star$ on $h_*(Ad(E))$ of degree $-m$:
$$\begin{CD}
h_r(Ad(E)^{-T_v}) \times h_s(Ad(E)^{-TM}) @>>> h_{r+s}(Ad(E)^{-T_v}) \\
@V{\cong}VV @VV{\cong}V \\
h_{r+g}(Ad(E)) \times h_{s+m}(Ad(E)) @>{\star}>> h_{r+s+g}(Ad(E))
\end{CD} $$
  We will now check that this is the same (up to sign) as an ``inverted'' string product.  Let $p_i: Ad(E)\times Ad(E) \to Ad(E)$ and $\pi _i :Ad(E) \times _M Ad(E) \to Ad(E)$ be the projections onto the $i$-th factor.  Take $a \in h_{r+g}(Ad(E))$ and $b \in h_{s+m}(Ad(E))$.  Then 
\begin{align*}
a \star b &= \pm {p_1}_* (\tilde{\Delta}^! D(a) \cap ([Ad(E)]\times b)) \\
&= \pm {p_1}_*\tilde{\Delta}_* (m^*D(a) \cap \tilde{\Delta}_!([Ad(E)]\times b) \\
&=\pm {p_1}_* \tilde{\Delta}_*((m\times id)^*(\pi_1^*D(a) \cup \pi_2^*D(b)) \cap (m \times id)_![Ad(E) \times _M Ad(E)]) \\
&= \pm {\pi_1}_* (m\times id)_! ((\pi_1^*D(a) \cup \pi_2^*D(b) ) \cap [Ad(E) \times _M Ad(E)]) \\
&= \pm {\pi_1}_* (m\times id)_! (\tilde{\Delta}_!(a \times b)) \\
&= \pm (\mbox{deg}(m \times id)){\pi_1}_*(\bar{m} \times id)_*(\tilde{\Delta}_!(a\times b)) \\
&= \pm {\bar{m}}_*(\tilde{\Delta}_!(a\times b)
\end{align*}
%
%
Here $\bar{m}: Ad(E)\times _M Ad(E) \to Ad(E)$ is defined by $\bar{m}(x,y) = m(x,y^{-1}).$  Notice that in general $\mbox{deg}(f) \in h_*(\ast)$, but since $(\mbox{deg}(m \times id))^2 = 1$, we know $\mbox{deg}(m \times id) = \pm 1$.  
We will not need to use the precise sign on the last line, just the fact that it depends only on the dimensions $r, s, g,$ and $m$.

The multiplication $\star$ is not associative in general, so does not define a Frobenius algebra structure with the comultiplication.  However, $\star$ does satisfy the property that 
\begin{equation} \label{bimod}
\begin{CD}
h_*(Ad(E)) \times h_*(Ad(E)) @>{\star}>> h_*(Ad(E)) \\
@V{\rho \times 1}VV @VV{\rho}V \\
h_*(Ad(E)) \times h_*(Ad(E)) \times h_*(Ad(E)) @>{1 \times \star}>> h_*(Ad(E))
\end{CD}
\end{equation}
commutes. 
  This follows from the homotopy-commutativity of the diagram below:
$$\begin{CD}
Ad(E)^{-T_v} \wedge Ad(E)^{-TM} @>{\rho \wedge 1}>> Ad(E)^{-T_v} \wedge Ad(E)^{-T_v} \wedge Ad(E)^{-TM} \\
@V{\rho \wedge 1}VV @VV{1\wedge \rho \wedge 1}V \\
Ad(E)^{-T_v}\wedge Ad(E)^{-T_v}\wedge Ad(E)^{-TM} @>{\rho \wedge 1 \wedge 1}>> Ad(E)^{-T_v} \wedge Ad(E)^{-T_v}\wedge  Ad(E)^{-T_v} \wedge Ad(E)^{-TM} \\
@V{1 \wedge \eta}VV @VV{1 \wedge 1 \wedge \eta}V \\
Ad(E)^{-T_v}\wedge S @>{\rho \wedge 1}>> Ad(E)^{-T_v}\wedge Ad(E)^{-T_v} \wedge S
\end{CD}$$

The top square commutes by the homotopy co-associativity of the comultiplication, which follows from the homotopy associativity of the multiplication on $Ad(E)^{-TM}$.  Our goal is to show that the diagram (\ref{bimod}) commutes with $\star$ replaced by $\mu$.  The commutativity of (\ref{bimod}) is the statement that for any $a, b \in h_*(Ad(E))$,
$$\tilde{\Delta}_*m_!\bar{m}_*\tilde{\Delta}_!(a\times b) = (1\times \bar{m}_*\tilde{\Delta}_! )(\tilde{\Delta}_*m_!(a) \times b).$$
Let $i$ denote the automorphism of $Ad(E)$ that sends $x$ to $x^{-1}.$  Replacing $b$ by $i_*(b)$ and observing that $i_* = (\mbox{deg}(i))i_!$ yields on the left side:
\begin{align*}
\tilde{\Delta}_*m_!\bar{m}_*\tilde{\Delta}_!(a\times i_*b) &=
\tilde{\Delta}_*m_!\bar{m}_*\tilde{\Delta}_!(1 \times i)_!(a\times b) \\
&= \tilde{\Delta}_*m_!\bar{m}_*(1 \times i)_!\tilde{\Delta}_!(a\times b) \\
&= (\mbox{deg}(i)) \tilde{\Delta}_*m_!\bar{m}_*(1 \times i)_* \tilde{\Delta}_!(a\times b) \\
&= (\mbox{deg}(i)) \tilde{\Delta}_*m_!(m)_*\tilde{\Delta}_!(a\times b)
\end{align*}

The same argument applied to the right hand side shows
$$(1\times \bar{m}_*\tilde{\Delta}_! )(\tilde{\Delta}_*m_!(a) \times i_*b) = (\mbox{deg}(i))(1\times m_*\tilde{\Delta}_! )(\tilde{\Delta}_*m_!(a) \times b)$$

Hence $$  \tilde{\Delta}_*m_!m_*\tilde{\Delta}_!(a\times b) = (1\times m_*\tilde{\Delta}_! )(\tilde{\Delta}_*m_!(a) \times b) .$$
In other words, the comultiplication is a map of right $h_*(Ad(E))$-modules.  Repeating the argument using the left $Ad(E)^{-TM}$-module structure on $Ad(E)^{-T_v}$ shows that the comultiplication is also a map of left $h_*(Ad(E))$-modules.  Thus the Chas-Sullivan product and FHT coproduct give a Frobenius algebra structure on $h_*(Ad(E))$.
\end{proof}

 There is also a Frobenius algebra structure on $h^*(Ad(E))$ that is Poincar\'{e} dual to the one described on $h_*(Ad(E))$.

\section{Further Remarks}\label{section:further}
In \cite{isak}, Christensen and Isaksen describe an extension of the Spanier-Whitehead duality functor that carries spectra to pro-spectra, giving a Quillen equivalence between the stable model structure on symmetric spectra \cite{hss} and a model structure on pro-spectra where the cohomotopy groups detect weak equivalences.  Their functor rectifies the problem of the poor behavior of Spanier-Whitehead duality for infinite spectra.  The Spanier-Whitehead dual of an infinite spectrum is then the pro-spectrum created from the duals of its finite subcomplexes.  From this point of view, Theorem~\ref{dual} suggests that $LBG^{-TBG}$ is the Spanier-Whitehead dual of $Ad(EG)^{-T_v}$.  However, $LBG^{-TBG}$ is a pro-ring spectrum in a weak sense, meaning it is a pro-object in a homotopy category of spectra, each spectrum $Ad(EG_n)^{-TM_n}$ is a ring object in the homotopy category, and the maps between them are homotopy ring maps.  In particular, in defining $LBG^{-TBG}$, we have not specified an object in Christensen and Isaksen's category of pro-spectra.  Keeping this is mind, we have the following corollary to Theorem~\ref{dual}.

\begin{cor}\label{bigdual}
$LBG^{-TBG}$ is  equivalent to the Spanier-Whitehead dual of $Ad(EG)^{-T_v}$ in a weak sense.
\end{cor}

Here what we mean by ``in a weak sense" is that the weak equivalences of Theorem~\ref{dual} homotopy commute with the structure maps of the pro-spectra.

We end with a conjecture about how this structure may be rigidified.  It is possible to rigidify $LBG^{-TBG}$ to a pro-object in the category of symmetric or orthogonal spectra \cite{Gthesis}.  

\begin{conj}\label{strictbigdual}
There is a rigidification of $LBG^{-TBG}$ which is weakly equivalent in the model category of \cite{isak} to the Spanier-Whitehead dual of $Ad(EG)^{-T_v}$.
\end{conj}

Notice that $Ad(EG)^{-T_v}$ still has a co-ring structure up to homotopy, defined just as in (\ref{coprod}).  Applying $K$-theory and Thom isomorphisms to this coproduct structure yields the product on $K^*(Ad(EG)) \cong K_G^*(G)\hat{}\,$ that is the completion of the fusion product on $K_G^*(G)$, as in (\ref{completion}).  Thus the moral of Corollary~\ref{bigdual} is that the string topology of $BG$ is Spanier-Whitehead dual to the completion of the fusion product on $K_G^*(G)$.

\bibliography{biblio}

\providecommand{\bysame}{\leavevmode\hbox to3em{\hrulefill}\thinspace}
\providecommand{\MR}{\relax\ifhmode\unskip\space\fi MR }
\providecommand{\MRhref}[2]{%
  \href{http://www.ams.org/mathscinet-getitem?mr=#1}{#2}
}
\providecommand{\href}[2]{#2}
\begin{thebibliography}{FHT03}

\bibitem[Abr96]{abrams}
L.~Abrams, \emph{Two dimensional topological quantum field theories and
  {F}robenius algebras}, J. Knot Theory Ramifications \textbf{5} (1996),
  569--587.

\bibitem[AS69]{as}
M.~F. Atiyah and G.~B. Segal, \emph{Equivariant {K}-theory and completion}, J.
  Differential Geom. \textbf{3} (1969), 1--18.

\bibitem[AS04]{astwk}
Atiyah-Segal, \emph{Twisted {K}-theory}, preprint: math.KT/0407054 (2004).

\bibitem[Ati61]{a}
M.~F. Atiyah, \emph{Thom complexes}, Proc. London Math. Soc. \textbf{11}
  (1961), no.~3, 291--310.

\bibitem[BG75]{bg}
J.~C. Becker and D.~H. Gottlieb, \emph{The transfer map and fiber bundles},
  Topology \textbf{14} (1975), 1--12.

\bibitem[CG04]{cg}
R.~L. Cohen and V.~Godin, \emph{A polarized view of string topology}, Topology,
  Geometry, and Quantum Field Theory, London Math. Soc. lecture notes, vol.
  308, 2004, pp.~127--154.

\bibitem[CI04]{isak}
J.~D. Christensen and D.~C. Isaksen, \emph{Duality and pro-spectra}, Algebr.
  Geom. Topol. \textbf{4} (2004), 781--812.

\bibitem[CJ02]{cj}
R.~L. Cohen and J.~D.~S. Jones, \emph{A homotopy theoretic realization of
  string topology}, Math. Ann. \textbf{324 no. 4} (2002), 773--798.

\bibitem[CS99]{cs}
M.~Chas and D.~Sullivan, \emph{String topology}, preprint: math.GT/9911159
  (1999).

\bibitem[FHT03]{fht3}
D.~S. Freed, M.~J. Hopkins, and C.~Teleman, \emph{Loop groups and twisted
  {K}-theory {III}}, preprint: math.AT/0312155 (2003).

\bibitem[Fre01]{freedv}
D.~S. Freed, \emph{The {V}erlinde algebra is twisted equivariant {K}-theory},
  Turk. J. Math \textbf{25} (2001), 159--168.

\bibitem[Gru]{Gthesis}
K.~Gruher, Ph.D. thesis, Stanford University, in preparation.

\bibitem[GS06]{gs}
K.~Gruher and P.~Salvatore, \emph{Generalized string topology operations},
  preprint: math.AT/0602210 (2006).

\bibitem[GW06]{gw}
K.~Gruher and C.~Westerland, \emph{String topology and homotopy fixed point
  prospectra}, in preparation (2006).

\bibitem[HSS00]{hss}
M.~Hovey, B.~Shipley, and J.~Smith, \emph{Symmetric spectra}, J. Amer. Math.
  Soc. \textbf{13} (2000), no.~1, 149--208.

\end{thebibliography}

\end{document}